\documentclass[12pt]{article}

\usepackage{doublespace,amstex}

\newtheorem{theorem}{Theorem}[section]
\newtheorem{lemma}[theorem]{Lemma}
\newtheorem{proposition}[theorem]{Proposition}
\newcommand{\gl}{\mathit{gl}}
\newcommand{\GL}{\mathit{GL}}
\newcommand{\GAP}{\mathit{GAP}}
\newtheorem{definition}[theorem]{Definition}
\newtheorem{corollary}[theorem]{Corollary}
\newtheorem{example}[theorem]{Example}
\newtheorem{problem}[theorem]{Problem}

\newcommand{\bproof}{\noindent{\bf Proof: }}
\newcommand{\eproof}{\hfill $\Box$\\}

\newcommand{\w}{\omega}

\newcommand{\e}{\varepsilon}

\newcommand{\cA}{{\cal A}}
\newcommand{\cS}{{\cal S}}
\newcommand{\cB}{{\cal B}}
\newcommand{\bN}{{\Bbb N}}
\newcommand{\bR}{{\Bbb R}}

\begin{document}

\title{A Gaussian Average Property of Banach Spaces}
 
\author{P.G.\ Casazza\thanks{This research was carried out during the first named author's visit to
Odense University in the Spring of 1995. It was partially supported by
the American National Science Foundation, grant no.\ NSF DMS-9201357 and by the Danish Natural Science Research Council, grant no.\ 9401598.} \and N.J.\ Nielsen
}

\date{ }
\maketitle
 
\renewcommand{\theequation}{\arabic{section}.\arabic{equation}}

\begin{abstract}

In this paper we investigate a Gaussian average property of Banach
spaces. This property is weaker than the Gordon Lewis property but
closely related to this and other unconditional structures. It is also
shown that this property implies that certain Hilbert space valued
operators can be extended.
\end{abstract}

\section*{Introduction}

In this paper we introduce a Gaussian average property, abbreviated
$\GAP$. A Banach space $X$ is said to have $\GAP$ if there is a constant $K$
so that $\ell(T)\le K\pi_1(T^*)$ for every finite rank operator from
$\ell_2$ to $X$. Here $\ell(T)$ denotes the $\ell$-norm defined by Linde
and Pietsch \cite{LIP}, see also N.\
Tomczak-Jaegermann \cite{NT}.

We investigate this property in detail and establish that a large class
of Banach spaces have it. It turns out that every Banach space, which is
either of type 2 or is isomorphic to a subspace of a Banach lattice of
finite cotype, has $\GAP$ and so does a Banach space of finite cotype which
has the Gordon-Lewis property $\GL_2$ with respect to  Hilbert spaces.

Though $\GAP$ is weaker than $\GL_2$, they are closely related, and since
$\GAP$ is somewhat easier to work with, this enables us to obtain some
results on $\GL_2$ by investigating $\GAP$. We prove e.g.\ that $\GAP$ and
$\GL_2$ are equivalent properties for cotype 2 spaces and that a
$K$-convex Banach space $X$ has $\GL_2$ if and only if both $X$ and $X^*$
have $\GAP$. It also turns out that if a space $X$ is of finite cotype and
$X^*$ has $\GAP$, then $X$ is $K$-convex.

We also prove that $\GAP$ gives rise to some extension theorems of
operators with range in a Hilbert space. We prove e.g.\ that if $X$ has
$\GAP$, then every operator from a subspace of $X$ into a Hilbert space,
which factors through $L_1$, extends to an $L_1$-factorable operator
defined on $X$. Further, if the dual of a subspace $E$ of a finite
cotype Banach space $X$ has $\GAP$, then every absolutely summing operator
from $E$ to a Hilbert space extends to an absolutely summing operator
defined on $X$. If $X^*$ has $\GAP$ then the other direction is true for
all subspaces $E$ of $X$. This implies that if $X$ is a
Banach space of finite cotype with $\GL_2$ then a subspace $E$ has $\GL_2$ if and only if
every 1-summing operator from $E$ to a Hilbert space extends to a
1-summing operator defined on $X$.

We now wish to discuss the arrangement and contents of the paper in
greater detail.

In Section \ref{sec-1} we prove the major results on $\GAP$ mentioned
above. One of the main tools for obtaining these is the duality theorem
\ref{thm-1.6} which also relates $\GAP$ to $K$-convexity. We provide several examples
of Banach spaces with a reasonable structure which fail $\GAP$. At the end
of the section it is shown that the $\ell_2$-sum of a sequence of
Banach spaces with uniformly bounded $\GAP$-constants (respectively
uniformly bounded $\GL_2$-constants) has $\GAP$ (respectively $\GL_2$). This
is obtained from an inequality for $p$-summing operators defined on an
$\ell_2$-sum of a sequence of Banach spaces with values in a Hilbert
space (Theorem \ref{thm-1.17}), which turns out to have applications also outside
the scope of this paper.

Section \ref{sec-2} is devoted to the extension theorems mentioned above.

\setcounter{section}{-1}
\section{Notation and Preliminaries}
\label{sec-0}

In this paper we shall use the notation and terminology commonly used in
Banach space theory, as it appears in \cite{L-T1}, \cite{L-T2} and
\cite{NT}.

If $X$ and $Y$ are Banach spaces, $B(X,Y)$ ($B(X)=B(X,X)$) denotes the
space of bounded linear operators from $X$ to $Y$. Further, if $1\le
p<\infty$ we let $\Pi_p(X,Y)$ denote the space of $p$-summing operators
from $X$ to $Y$ equipped with the $p$-summing norm $\pi_p$. We recall
that an operator $T\in B(X,Y)$ is said to factor through $L_p$ if it
admits a factorization $T=BA$, where $A\in B(X,L_p(\nu))$ and $B\in
B(L_p(\nu),Y)$ for some measure $\nu$ and we denote the space of all
operators which factor through $L_p$ by $\Gamma_p(X,Y)$. If
$T\in\Gamma_p(X,Y)$ then we define
\[
\gamma_p(T) = \inf\{\|A\| \|B\| \mid T = BA,\quad \mbox{$A$ and $B$ as
above}\}.
\]
$\gamma_p$ is a norm on $\Gamma_p(X,Y)$ turning it into a Banach space.

Throughout the paper we shall identify the tensor product $X\otimes Y$
with the space of $\w^*$-continuous finite rank operators from $X^*$ to
$Y$ in the canonical manner.

We let $(g_n)$ denote a sequence of independent standard Gaussian
variables on a fixed probability space $(\Omega,\cS,\mu)$ and we let
$G(X)$ denote the closure of $\left\{\sum_{j=1}^n g_j x_j\mid n\in\bN \
x_j\in X\right.$, $\left.1\le j\le n\right\}$ in $L_2(\mu,X)$. Further, we let $(e_n)$
denote the unit vector basis of $\ell_2$.

If $n\in\bN$ and $T\in B(\ell_2^n,X)$ then, following \cite{NT}, we
define the $\ell$-norm of $T$ by
\[
\ell(T) = \left(\int \left\|\sum_{j=1}^n g_j(t) Te_j\right\|^2
d\mu(t)\right)^{\frac12}.
\]

More generally, if $T\in B(\ell_2,X)$, we call $T$ an $\ell$-operator if
$\sum_{n=1}^\infty g_n Te_n$ converges in $L_2(\mu,X)$ and we put
\[
\ell(T) = \left(\int\left\|\sum_{n=1}^\infty g_n(t)Te_n\right\|^2
d\mu(t)\right)^{\frac12}.
\]

We also need some notation on operators with ranges in a Banach
lattice. Recall that if $E$ is a Banach space and $X$ is a Banach
lattice, then an operator $T\in B(E,X)$ is called order bounded (see
e.g.\ \cite{NJN1} and \cite{NJN2}), if there exists a $z\in X$, $z\ge 0$
so that
\begin{equation}
\label{eq-*}
|Tx|\le \|x\|z \quad\mbox{for all $x\in E$}
\end{equation}
and the order bounded norm $\|T\|_m$ is defined by
\begin{equation}
\label{eq-**}
\|T\|_m = \inf\{\|z\| \mid \mbox{$z$ can be used in (\ref{eq-*})} \}.
\end{equation}
$\cB(E,X)$ denotes the Banach space of all order bounded operators from
$E$ to $X$ equipped with the norm $\|\cdot\|_m$.

A Banach space $X$ is said to have the Gordon-Lewis property
(abbreviated $\GL$) \cite{G-L} if every 1-summing operator from $X$ to an
arbitrary Banach space $Y$ factors through $L_1$. It is readily verified
that $X$ has $\GL$, if and only if there is a constant $K$ so that
$\gamma_1(T)\le K\pi_1(T)$ for every Banach space $Y$ and every $T\in
X^*\otimes Y$. In that case $\gl(X)$ denotes the smallest constant $K$
with this property.

We shall say that $X$ has $\GL_2$ if it has the property above with
$Y=\ell_2$ and we define the constant $\gl_2(X)$ correspondingly. An
easy trace duality argument yields that $\GL$ and $\GL_2$ are self dual
properties and that $\gl(X)=\gl(X^*)$ when applicable. It is known
\cite{G-L} that every Banach space with local unconditional structure
has the Gordon-Lewis property.

We now present a few theorems which all follow from well-known results
and which do not appear in the literature in the form we are going to
use them.

The first proposition follows immediately from the contraction principle
for independent Gaussian variables, see e.g.\ \cite{P2}.

\begin{proposition}
\label{prop-0.1}
If $X$ is a Banach space and $(x_j)\subseteq X$ then for all $n\le m$
and all $1\le p<\infty$ we have
\[
\int\left\|\sum_{j=1}^n g_j(t)x_j\right\|^p d\mu(t) \le
\int\left\|\sum_{j=1}^m g_j(t)x_j\right\|^p d\mu(t).
\]
\end{proposition}

As a corollary to Proposition \ref{prop-0.1} we obtain

\begin{proposition}
\label{prop-0.2}
If $X$ is a Banach space and $f\in G(X)$ then for all $n\in\bN$ we have
\begin{equation}
\label{eq-i}
\left(\int\left\|\sum_{j=1}^n g_j(t)\int f\,
g_jd\mu\right\|^2d\mu(t)\right)^{\frac12} \le \|f\|_2,
\end{equation}
and the series $\sum_{j=1}^\infty g_j(\int f\, g_jd\mu)$ converges to
$f$ in $L_2(\mu,X)$.
\end{proposition}

\bproof
Let $\cA$ be the subspace of $G(X)$ consisting of all $f$ of the form
$f=\sum_{j=1}^m g_jx_j$ for some $m\in\bN$ and some sequence
$(x_j)\subseteq X$. For every $n\in\bN$ we define $P_n: \cA\to G(X)$ by
\begin{equation}
\label{eq-1}
P_nf = \sum_{j=1}^n g_j\left(\int fg_jd\mu\right) \ \mbox{for all
$f\in\cA$}.
\end{equation}
>From the previous proposition it follows that $P_n$ is a bounded linear
projection on $\cA$ with $\|P_n\|\le 1$. Hence it can be extended to a
norm 1 linear projection on $G(X)$ also denoted $P_n$. This gives
immediately (\ref{eq-i}) and since obviously $P_nf\to f$ for all $f\in\cA$
and $\|P_n\|\le 1$ for all $n\in\bN$ we get the same for all $f\in G(X)$.
\eproof

\begin{proposition}
\label{prop-0.3}
For every $1\le p<\infty$ there is a constant $K_p$ so that if $X$ is a
Banach space and $T\in B(\ell_2,X)$ is an $\ell$-operator then $T^*$ is
$p$-summing with
\begin{equation}
\label{eqi}
\pi_p(T^*)\le K_p \ell(T).
\end{equation}
If $T\in\ell_2\otimes X$, then $T$ is an $\ell$-operator and
\begin{equation}
\label{eqii}
\ell(T)\le K_p\pi_p(T).
\end{equation}
\end{proposition}

\bproof
Let $1\le p<\infty$. By a result of Kahane \cite{KA} there are constants
$a_p>0$ and $b_p>0$ so that
\begin{equation}
\label{eq1}
a_p\|f\|_2\le \|f\|_p \le b_p\|f\|_2 \quad\mbox{for all $f\in G(X)$}.
\end{equation}
To prove (\ref{eqi}) we let $T\in B(\ell_2,X)$ be an $\ell$-operator and
define
\begin{equation}
\label{eq2}
f=\sum_{n=1}^\infty g_nTe_n.
\end{equation}

If $I_p:$ $\ell_2\to L_p(\mu)$ denotes the operator defined by $I_p
e_n=g_n$ for all $n\in\bN$, then $I_p$ is an isomorphism and
\begin{equation}
\label{eq3}
(I_pT^*x^*)(t) = x^*(f(t)) \quad\mbox{for almost all $t\in\Omega$}.
\end{equation}

It follows from \cite{NJN1} that $I_pT^*$ is order bounded and therefore
$p$-summing with
\begin{equation}
\label{eq4}
\pi_p(T^*)\le a_p^{-1}\pi_p(I_pT^*)\le a_p^{-1}\|I_pT^*\|_m =
a_p^{-1}\|f\|_p \le a_p^{-1}b_p\ell(T).
\end{equation}
To prove (\ref{eqii}) we let $T\in\ell_2\otimes X$; hence there is a
$K\in\bN$, $(f_j)_{j=1}^k \subseteq\ell_2$ and $(x_j)_{j=1}^k\subseteq X$
with $T=\sum_{j=1}^k f_j\otimes x_j$. If $g=\sum_{j=1}^k (I_2f_j)x_j$,
then for all $n\in\bN$,
\begin{equation}
\label{eq5}
Te_n=\sum_{j=1}^k (e_n,f_j)x_j = \sum_{j=1}^k (g_n,I_2f_j)x_j = \int
g(t)g_n(t)d\mu(t),
\end{equation}
and therefore, by Proposition \ref{prop-0.2}
\begin{equation}
\label{eq6}
g=\sum_{n=1}^\infty g_nTe_n.
\end{equation}

This shows that $T$ is an $\ell$-operator and by \cite[Corollary
4.8]{NJN1} we obtain
\begin{align}
\label{eq7}
\|I_pT^*\|_m&\le\pi_p(T)\\
\intertext{and}
\label{eq8}
\ell(T)\le a_p^{-1}\|I_pT^*\|_m&\le a_p^{-1}\pi_p(T).
\end{align}
\eproof


\setcounter{equation}{0}

\section{The Gaussian Average Property and Related Topics}
\label{sec-1}

In this section we shall introduce our Gaussian average property and
prove our main results, which among other things relates this property
to the Gordon-Lewis property. We start with the following definition

\begin{definition}
\label{def-1.1}
Let $X$ be a Banach space. $X$ is said to have the Gaussian average
property ($\GAP$) if there is a constant $K$, so that for all
$T\in\ell_2\otimes X$ we have $\ell(T)\le K\pi_1(T^*)$.

$X$ is said to have property $(S_p)$ $1\le p<\infty$ if there is a
constant $K$ so that if $T\in B(\ell_2,X)$ with
$T^*\in\Pi_1(X^*,\ell_2)$, then $T\in\Pi_p(\ell_2,X)$ with $\pi_p(T)\le
K\pi_1(T^*)$. We shall say that $X$ has $(S)$, if it has $(S_p)$ for some
$p$, $1\le p<\infty$.
\end{definition}

Recall that a Banach space $Y$ is called a Grothendieck space
(abbreviated GT) \cite{P2} if $B(Y,\ell_2)=\Pi_1(Y,\ell_2)$. It follows
from Grothendieck's inequality that every ${\cal L}_1$-space is a GT
space. We make the following observation:

\begin{proposition}
\label{prop-1.2a}
If $X$ is a Banach space so that $X^*$ is a GT-space then $X$ does not
have $\GAP$. In particular, $L_\infty$ does not have $\GAP$.
\end{proposition}

\bproof
Let $K$ be the GT-constant of $X$ and let $n\in\bN$ be given. By
Dvoretzky's theorem \cite{L-T1} there is an isomorphism $T: \ell_2^n\to
X$ so that $\|T\|\le 2$ and $\|T^{-1}\|=1$. Clearly $\pi_1(T^*)\le
K\|T\|\le 2K$ and $\frac12 \sqrt{n}\le \ell(T)\le 2\sqrt{n}$, which
shows that $X$ does not have $\GAP$.
\eproof   

It follows easily from the results of the previous section that if $X$
as $\GAP$, then the $\ell$-norm of an operator $T\in B(\ell_2,X)$ is
equivalent to the 1-summing norm of the adjoint. If $X$ has $(S_p)$ then
it follows that the $p$-summing norm of an operator
$T\in\Pi_p(\ell_2,X)$ is equivalent to the 1-summing norm of the
adjoint.

It is readily seen that both $\GAP$ and $(S)$ are hereditary properties and
from the principle of local reflexivity it is easily seen that $X$ has
$\GAP$, respectively $(S)$, if and only if $X^{**}$ has $\GAP$, respectively
$(S)$. Furthermore we have

\begin{theorem}
\label{thm-1.2}
Let $X$ be a Banach space. Then the following statements hold:
\begin{itemize}
\item[(i)] If $X$ has $(S)$, then it has $\GAP$.
\item[(ii)] If $X$ has $(S_p)$, then it is of cotype $\max(2,p)$.
\item[(iii)] If $X$ has $\GAP$, then it is of finite cotype.
\item[(iv)] If $X$ is of finite cotype and has $\GL_2$, then $X$ has
$(S)$ and hence also $\GAP$.
\end{itemize}
\end{theorem}

\bproof
(i) and (ii): Let $X$ have $(S_p)$ with constant $K$ for some $p$, $1\le
p<\infty$ and put $q=\max(p,2)$. It follows from Proposition
\ref{prop-0.3} that for every $T\in \ell_2\otimes X$ we have:
\begin{equation}
\label{eq(1)}
\pi_{q,2}(T) \le \pi_p(T)\le K\pi_1(T^*)\le KK_1\ell(T)\le
K_pKK_1\pi_p(T)\le K^2K_pK_1\pi_1(T^*).
\end{equation}
>From (\ref{eq(1)}) we obtain directly that $X$ has $\GAP$. Furthermore,
together with \cite[Theorem 12.2]{NT}, (\ref{eq(1)}) gives that $X$ has cotype $q$.

(iii): Assume that $X$ has $\GAP$. If $X$ is not of finite cotype it
contains $\ell^n_\infty$ uniformly \cite{M-P} and since $\GAP$ is hereditary this
implies that $\ell_\infty$ has $\GAP$, which is a contradiction.

(iv): Let $X$ be a Banach space of cotype $q$ with $\GL_2$ and let
$p>q$. By self-duality $X^*$ has $\GL_2$ as well and if $T\in
B(\ell_2,X)$ with $T^*\in\Pi_1(X^*,\ell_2)$ then
$T\in\Gamma_\infty(\ell_2,X^{**})$ 
and hence by \cite{M-P} 
$T\in\Pi_p(\ell_2,X)$. If $q=2$, we can actually take $p=2$ as well.
\eproof

The next theorem describes some classes of Banach spaces which have $\GAP$.

\begin{theorem}
\label{thm-1.3}
Let $X$ be a Banach space. 
\begin{itemize}
\item[(i)] If $X$ is of cotype 2 then $X$ has $\GAP$ if and only if
it has $\GL_2$.
\item[(ii)] If $X$ is of type 2 then it has $\GAP$.
\item[(iii)] If $X$ is a subspace of a Banach lattice of finite cotype,
then $X$ has $(S)$ and hence $\GAP$.
\end{itemize}
\end{theorem}

\bproof
(i): If $X$ is of cotype 2 it follows from \cite[Theorem 12.2]{NT} that
there is a constant $K$ so that
\begin{equation}
\label{eq-1.1}
\pi_2(T) \le K\ell(T) \quad\mbox{for all $T\in\ell_2\otimes X$}.
\end{equation}

If $X$ has $\GAP$ with constant $C$ then it follows from (\ref{eq-1.1})
that for all $T\in\ell_2\otimes X$ we have
\begin{equation}
\label{eq-1.2}
\gamma_1(T^*) = \gamma_\infty(T)\le\pi_2(T)\le K\ell(T)\le
KC\pi_1(T^*).
\end{equation}
This shows that $X^*$ and hence $X$ has $\GL_2$.

The other direction follows from Theorem \ref{thm-1.2}.

(ii): Let $X$ be of type 2 with constant $K$ and let $T\in\ell_2\otimes
X$. Again, by \cite[Theorem 12.2]{NT}, we get:
\begin{equation}
\label{eq-1.3}
\ell(T)\le K\pi_2(T^*)\le K\pi_1(T^*),
\end{equation}
which shows that $X$ has $\GAP$.

(iii): Let $X$ be a subspace of a Banach lattice $Z$ of finite
cotype. Hence, by \cite{L-T2}, $Z$ is $q$-concave for some $q$, $1\le
q<\infty$ with constant say $K$. If $T\in\ell_2\otimes X$ and $I:$ $X\to
Z$ denotes the identity operator, then it follows from \cite[Proposition 4.9]{NJN1} that
\begin{equation}
\label{eq-1.4}
\|IT\|_m\le \pi_1(T^*I^*)\le\pi_1(T^*).
\end{equation}
Since $T$ is of finite rank it follows from \cite[Theorem 2.9]{NJN1} that
there exists a compact Hausdorff space $S$ and operators $A\in
B(\ell_2,C(S))$, $B\in B(C(S),Z)$ so that $\|A\|=1$, $B\ge 0$,
$\|B\|=\|IT\|_m$ and $IT=BA$. Since $B\ge 0$ and $Z$ is $q$-concave, $B$
is $q$-summing with $\pi_q(B)\le K\|B\|$ (\cite{L-T2}). Hence $T$ is
$q$-summing as well with
\begin{equation}
\label{eq-1.5}
\pi_q(T)\le \|A\|\pi_q(B)\le K\|T\|_m\le K\pi_1(T^*).
\end{equation}
This shows that $X$ has $(S_q)$.
\eproof

Since $\GAP$ is a hereditary property, Theorem \ref{thm-1.3} gives the
following corollary:

\begin{corollary}
\label{cor-1.4}
If $X$ of cotype 2 has $\GL_2$ then so does every subspace. In
particular, if $X$ is a Banach lattice of cotype 2, then every subspace
has $\GL_2$.
\end{corollary}

Corollary \ref{cor-1.4} can of course also easily be deduced from the
fact that if $X$ is of cotype 2 then $\Pi_1(X,L_2)=\Pi_2(X,L_2)$ and the
fact that 2-summing operators extend to 2-summing operators.

The cotype 2 situation is not the only one where $\GAP$ and $\GL_2$
coincide. We shall return to this after we have proved an important
duality theorem. First we need:

\begin{proposition}
\label{prop-1.5}
If $X$ is a Banach space of finite cotype then there is a constant $K\ge
0$ so that
\begin{equation}
\label{eq-prop1.5}
\ell(T)\le K\gamma_\infty(T) \quad\mbox{for all $T\in\ell_2\otimes X$}.
\end{equation}
\end{proposition}

\bproof
Let $X$ be of cotype $r$ and let $q>r$. From \cite{M-P} it is easily
derived that there is a constant $K$ so that
\begin{equation}
\label{eq-1.6}
\pi_q(T)\le K\gamma_\infty(T)\quad\mbox{for all $T\in\ell_2\otimes X$}.
\end{equation}
(\ref{eq-prop1.5}) now follows by combining (\ref{eq-1.6}) with Theorem 
\ref{thm-1.3}.
\eproof

We are now able to prove the following duality theorem

\begin{theorem}
\label{thm-1.6}
If $X$ is a Banach space then the following conditions are equivalent:
\begin{itemize}
\item[(i)] $X$ is $K$-convex and there is a constant $K\ge 0$ so that
$K^{-1}\gamma_\infty(T)\le \ell(T)\le K\gamma_\infty(T)$ for
all $T\in\ell_2\otimes X$.
\item[(ii)] $X^*$ has $\GAP$ and $X$ is of finite cotype.
\end{itemize}
\end{theorem}

\bproof
(i) $\Rightarrow$ (ii). Assume that (i) holds and let $C$ denote the
$K$-convexity constant of $X$ (for the definition of $K$-convexity we
refer to \cite{P2}).

If $S\in\ell_2\otimes X^*$ we get
\begin{eqnarray}
\label{eq-1.7}
\ell(S) &\le & C \sup\{|Tr(T^*S)| \mid T\in\ell_2\otimes X, \ \ell(T)\le 1\}
\nonumber \\
&\le & KC \sup\{|Tr(T^*S)|\mid T\in \ell_2\otimes X, \ \gamma_\infty(T)\le
1\} = KC\pi_1(S^*),
\end{eqnarray}
which shows that $X^*$ has $\GAP$. Clearly $X$ is of finite cotype.

(ii) $\Rightarrow$ (i). Since $X$ is of finite cotype it follows from
Proposition \ref{prop-1.5} that there is a constant $C_1$ so that
$\ell(T)\le C_1\gamma_\infty(T)$ for all $T\in \ell_2\otimes X$. If
$C_2$ denotes the $\GAP$-constant of $X^*$ we get for every $T\in
\ell_2\otimes X$
\begin{eqnarray}
\label{eq-1.8}
\gamma_\infty(T) &=& \sup\{|Tr(S^*T)| \mid S\in\ell_2\otimes
X^*,\pi_1(S^*)\le 1\} \nonumber \\
&\le & C_2\sup\{|Tr(S^*T)|\mid S\in\ell_2\otimes X^*,\ell(S)\le 1\} \nonumber \\
&\le & C_2\ell(T)\le C_1C_2\gamma_\infty(T).
\end{eqnarray}
This shows that the third and fourth entries in (\ref{eq-1.8}) are
equivalent, which clearly (see \cite{NT}) implies that $X$ is
$K$-convex. In addition (\ref{eq-1.8}) shows that $\gamma_\infty(T)\le
C_2\ell(T)$ for all $T\in\ell_2\otimes X$. Hence we have proved that
(ii) $\Rightarrow$ (i).
\eproof

Since $X$ has $\GAP$ if and only if $X^{**}$ has $\GAP$, as noted just after
Definition \ref{def-1.1}, it follows that the roles of $X$ and $X^*$ can
be interchanged in Theorem \ref{thm-1.6}.

Theorem \ref{thm-1.6} has several corollaries

\begin{corollary}
\label{cor-1.7}
If $X$ has $\GAP$ and $X^*$ is of finite cotype then $X$ is $K$-convex.
\end{corollary}

The next corollary we formulate as a theorem

\begin{theorem}
\label{thm-1.8}
Let $X$ be a Banach space. The following statements are equivalent:
\begin{itemize}
\item[(i)] $X$ has $\GL_2$ and both $X$ and $X^*$ are of finite cotype.
\item[(ii)] $X$ and $X^*$ have $\GAP$.
\end{itemize}
Under these circumstances $X$ is $K$-convex.
\end{theorem}

\bproof
(i) $\Rightarrow$ (ii). Since $\GL_2$ is a self dual property it follows
that both $X$ and $X^*$ have $\GAP$.

(ii) $\Rightarrow$ (i). Assume that (ii) holds. It follows from Theorem
\ref{thm-1.2} that both $X$ and $X^*$ are of finite cotype.

Since $X$ has $\GAP$ it follows from Theorem \ref{thm-1.6} that there is a
constant $K\ge 0$ so that for all $T\in\ell_2\otimes X^*$ we have
\begin{equation}
\label{eq-1.9}
\gamma_\infty(T)\le K\ell(T).
\end{equation}
If $C$ denotes the $\GAP$-constant of $X^*$ we get from (\ref{eq-1.9}) that
if $S\in X\otimes\ell_2$, then
\begin{equation}
\label{eq-1.10}
\gamma_1(S) = \gamma_\infty(S^*) \le K\ell(S^*)\le KC\pi_1(S)
\end{equation}
which shows that $X$ has $\GL_2$.
\eproof

It is well known that if $X$ is of cotype 2 then
$B(L_\infty,X)=\Pi_2(L_\infty,X)$ or equivalently
$\Pi_1(X,\ell_2)=\Pi_2(X,\ell_2)$ and it is an open question whether the
converse implication holds. Pisier \cite{P1} showed that this is the
case if $X$ has $\GL_2$. Here we prove a similar result using $\GAP$.

\begin{theorem}
\label{thm-1.10}
Let $X$ be a Banach space and $1< p\le 2$,
$\frac{1}{p}+\frac{1}{q}=1$. If $X$ has $\GAP$, then
$B(L_\infty,X^*)=\Pi_q(L_\infty,X^*)$ if and only if $X$ is of type
$p$-stable.

In particular, $X$ is of type 2 if and only if it has $\GAP$ and
$\Pi_1(X^*,\ell_2)=\Pi_2(X^*,\ell_2)$.
\end{theorem}

\bproof
If $X$ is of type $p$-stable then it follows from \cite{M-P} that
$B(\ell_\infty,X^*)=\Pi_q(\ell_\infty,X^*)$. Assume next that $X$ has
$\GAP$ with constant $M$ and that $B(L_\infty,X^*)=\Pi_q(L_\infty,X^*)$
with $K$-equivalence between the norms, hence also
$\Pi_1(X^*,\ell_2)=\Pi_p(X^*,\ell_2)$ with $K$-equivalence between the
norms.

If $T=\sum_{j=1}^k e_j\otimes x_j\in\ell_2\otimes X$, then
\begin{equation}
\label{eq-1.11}
\pi_p(T^*) \le \left(\sum_{j=1}^k \|x_j\|^p\right)^{\frac{1}{p}}
\sup\left\{\left(\sum_{j=1}^k|(z,e_j)|^q\right)^{\frac{1}{q}} \mid
z\in\ell_2 \ \|z\|_2\le 1\right\}\le \left(\sum_{j=1}^k
\|x_j\|^p\right)^{\frac{1}{p}}
\end{equation}
and therefore
\begin{equation}
\label{eq-1.12}
\left(\int\left\|\sum_{j=1}^k g_j(t)x_j\right\|^2 d\mu(t)\right)^{\frac12}
= \ell(T) \le M\pi_1(T^*)\le MK\pi_p(T^*)\le \left(\sum_{j=1}^k
\|x_j\|^p\right)^{\frac{1}{p}}.
\end{equation}
which shows that $X$ is of type $p$.

If $p=2$ we are done. If $p<2$ then by \cite{M-P} $\{p<2 \mid
\Pi_q(L_\infty,X^*)=B(L_\infty,X^*)\}$ is an open interval and therefore
$X$ is of type $p$-stable.
\eproof

Let us now look on a few examples:

\begin{example}
\label{ex-1.9}
Let $X$ be the space constructed by Pisier in \cite{P1}. Both $X$ and
$X^*$ are of cotype 2, but $X$ is not isomorphic to a Hilbert
space. Therefore $X$ is not $K$-convex and hence cannot have $\GAP$ nor
$\GL_2$ by Corollary \ref{cor-1.7}.
\end{example}

There exist $K$-convex Banach spaces of cotype 2 not having $\GAP$
(equivalently $\GL_2$), which the following example shows:

\begin{example}
\label{ex-1.10}
Let $2<p<\infty$. By Johnson \cite[Lemma 1]{J-S}, it
follows that there exists a subspace $X\subseteq L_p(0,1)$ which does
not have $\GL_2$. $X^*$ is $K$-convex but does not have $\GAP$ by
Theorem \ref{thm-1.8}. Hence it does not embed into a Banach
lattice of finite cotype. 
\end{example}

Similar arguments as in this example leads to

\begin{corollary}
\label{cor-1.12}
Let $X$ be a Banach space with $\GAP$. If $X^*$ embeds into a Banach
lattice of finite cotype, then $X$ has $\GL_2$.
\end{corollary}

>From the theorem of Johnson we can also conclude

\begin{corollary}
\label{cor-1.13}
Every Banach lattice of finite cotype which is not of weak cotype 2
contains a subspace $X$, so that $X^*$ does not embed into a Banach
lattice of finite cotype.
\end{corollary}

We can pose the following problem:

\begin{problem}
\label{prob-1.14}
Can the above mentioned theorem of Johnson be
strengthened. Specifically, is a Banach space of cotype 2, if all
subspaces have $\GL_2$?
\end{problem}

The convexified Tsirelson space $T^{(2)}$ (see \cite{P2}) is of type 2
and weak cotype 2, and one could try to investigate whether there is a
subspace $X$ of $T^{(2)}$ failing $\GL_2$. Hence $X^*$ will fail $\GAP$ and
therefore $X$ would be the first example of a weak Hilbert space, which
does not embed any Banach lattice of finite cotype.

In \cite{Pisierny} Pisier showed that the Schatten class $c_p$, $p\ne 2$
does not have lust, but his argument actually gives that it does not
have (S). Indeed, in his Theorem 1.1 he shows that a space with lust
has (S) (called (I) there) and an inspection of the proof of
Proposition 2.1 in the paper shows that it is enough for the conclusion
to assume that the space $E$ there has (S). This observation together
with his proof of Theorem 2.1 then shows that if $\lambda$ is a
unitarily invariant crossnorm on $\ell_2\otimes\ell_2$ then
$\ell_2\hat{\otimes}_\lambda \ell_2$ does not have (S) unless
$\lambda$ is equivalent to the Hilbert Schmidt norm. In particular $c_p$
does not have (S) for $p\ne 2$.

Combining this with our Theorems \ref{thm-1.3} and \ref{thm-1.8} we
obtain:

\begin{example}
\label{ex-1.15a}
For every $q$, $2<q<\infty$, $c_q$ has $\GAP$ but not (S). $c_p$ does not
have $\GAP$ for $1\le p<2$.
\end{example}

The following condition is stronger than (S).

\begin{definition}
\label{def-1.16}
A Banach space $X$ is said to have $(I)$, if there is a $p$, $1\le
p<\infty$ and a constant $K$ so that
\[
i_p(T)\le K\pi_1(T^*)\quad \mbox{for all $T\in\ell_2\otimes X$}
\]
where $i_p$ denotes the $p$-integral norm \cite{NT}.
\end{definition}

Condition $(I)$ is equivalent to $X$ being of finite cotype and having
$\GL_2$. Indeed, if $X$ has $(I)$, then it has $(S)$ and is of finite
cotype. $(I)$ immediately implies that
$\Pi_1(X^*,\ell_2)\subseteq\Gamma_1(X^*,\ell_2)$ and therefore $X^*$ and
hence $X$ has $\GL_2$. On the other hand, if $X$ is of finite cotype and
has $\GL_2$, an inspection of the proof of Theorem \ref{thm-1.2}, (iv)
shows that in fact $X$ has $(I)$ (use that
$I_p(L_\infty,X)=\Pi_p(L_\infty,X)$ together with the principle of
local reflexivity).

This equivalence was also established by Junge \cite{J}.

We now wish to show that $\GAP$ is closed under the formation of
$\ell_2$-sums of Banach spaces. For this we need the following theorem,
which turns out to have some importance in itself.

\begin{theorem}
\label{thm-1.17}
Let $(X_n)$ be a sequence of Banach spaces and put $X=\left(\sum_{n=1}^\infty
X_n\right)_2$. If $Y$ is another Banach space, $1\le p<\infty$ and $T\in
\Pi_p(X,Y)$ with $T_n=T_{|X_n}$, then
\begin{eqnarray}
\left(\sum_{n=1}^\infty \pi_p(T_n)^2\right)^{\frac12} &\le &
\pi_p(T)\quad\mbox{for $1\le p\le 2$} \label{eq-1.13}
\\
\left(\sum_{n=1}^\infty \pi_p(T_n)^p\right)^{\frac{1}{p}} &\le &
\pi_p(T) \quad\mbox{for $2\le p<\infty$}.
\label{eq-1.131}
\end{eqnarray}
If $Y=\ell_2$ then (\ref{eq-1.13}) holds for all $p$, $1\le p<\infty$.
\end{theorem}

\bproof
Let $\e>0$ be given arbitrarily. For every $n\in\bN$ we can find a
finite set $\sigma_n\subseteq\bN$ and $\{x_i(n)\mid
i\in\sigma_n\}\subseteq X_n$ so that
\begin{equation}
\label{eq-1.14}
\pi_p(T_n)^p\le \sum_{i\in\sigma_n} \|Tx_i(n)\|^p +\e 2^{-n},
\end{equation}
\begin{equation}
\label{eq-1.15}
\sup\left\{\sum_{i\in\sigma_n} |x^*(x_i(n)|^p \mid x^*\in X^*_n,
\|x^*\|\le 1\right\}\le 1.
\end{equation}
For every sequence $(\alpha_n)\subseteq\bR_+\cup\{0\}$  we obtain from (\ref{eq-1.14}) and
(\ref{eq-1.15})
\begin{eqnarray}
\label{eq-1.16}
\lefteqn{\sum_{n=1}^\infty \alpha_n\pi_p(T_n)^p = \sum_{n=1}^\infty
\sum_{i\in\sigma_n} \|T(\alpha_n^{1/p}x_i(n))\|^p+\e}  \nonumber \\
&&\le 
 \pi_p(T)^p \sup\left\{\sum_{n=1}^\infty \sum_{i\in\sigma_n} |\langle
x^*(n),\alpha_n^{1/p}x_i(n)\rangle |^p \mid x^*(n)\in X_n^*,
\sum_{n=1}^\infty \|x^*(n)\|^2 \le 1\right\} +\e \nonumber \\
&&\le  \pi_p(T)^p \sup\left\{\sum_{n=1}^\infty \|x^*(n)\|^p \alpha_n
\sum_{i\in\sigma_n} |\langle \frac{x^*(n)}{\|x^*(n)\|},x_i(n)\rangle|^p
|x^*(n)\in X^*_n, \sum_{n=1}^\infty \|x^*(n)\|^2\le 1\right\}+\e 
\nonumber \\
&& \le \pi_p(T)^p \sup\left\{\sum_{n=1}^\infty \|x^*(n)\|^p \alpha_n\mid
x^*_n\in X^*_n, \sum_{n=1}^\infty \|x^*_n\|^2\le 1\right\} +\e 
\end{eqnarray}
If $1\le p\le 2$ we take the supremum in (\ref{eq-1.16}) over all
sequences $(\alpha_n)$ considered with $\sum_{n=1}^\infty
\alpha_n^{2/(2-p)}=1$ and let $\e\to 0$ to obtain (\ref{eq-1.13}).

For $2\le p<\infty$ we put $\alpha_n=1$ for all $n\in\bN$ in
(\ref{eq-1.16}) to obtain (\ref{eq-1.131}).

Since $\Pi_p(Z,\ell_2)=\Pi_2(Z,\ell_2)$ for every Banach space $Z$ and
every $2\le p<\infty$ (this follows easily from Maurey's extension
theorem \cite{M} and the formula $B(L_\infty,L_2)=\Pi_2(L_\infty,L_2)$)
the statement for $Y=\ell_2$ follows from the above.

\eproof

This enables us to prove

\begin{theorem}
\label{thm-1.18}
Let $(X_n)$ be a sequence of Banach spaces, which all have $\GAP$ with
uniformly bounded constants, then $X=\left(\sum_{n=1}^\infty
X_n\right)_2$ has $\GAP$.
\end{theorem}

\bproof
For every $n\in\bN$ we let $P_n$ denote the canonical projection of $X$
onto $X_n$. If $x_1,x_2,\dots,x_k\in X$ then it follows immediately from
the definition of the norm in $X$ that
\begin{equation}
\label{eq-1.17}
\int\left\|\sum_{i=1}^k g_1(t)x_i\right\|^2 d\mu(t) = \sum_{n=1}^\infty
\int \left\|\sum_{i=1}^k g_i(t)P_nx_i\right\|^2 d\mu(t).
\end{equation}
This gives that if $T\in B(\ell_2,X)$ is an $\ell$-operator, then
\begin{equation}
\label{eq-1.18}
\ell(T) = \left(\sum_{n=1}^\infty \ell(P_nT)^2\right)^{\frac12}.
\end{equation}
Let $K\ge 0$ be a constant so that for all $n\in\bN$,
\begin{equation}
\label{eq-1.19}
\ell(S)\le K\pi_1(S^*) \quad\mbox{for all $S\in\ell_2\otimes X$}.
\end{equation}

If now $T\in\ell_2\otimes X$, then by Theorem \ref{thm-1.17} with $p=1$,
(\ref{eq-1.18}) and (\ref{eq-1.19}), we obtain
\begin{equation}
\label{eq-1.20}
\ell(T) = \left(\sum_{n=1}^\infty \ell(P_nT)^2\right)^{\frac12} \le
K\left(\sum_{n=1}^\infty \pi_1(T^*P^*_n)^2\right)^{\frac12} \le
K\pi_1(T^*).
\end{equation}
This shows that $X$ has $\GAP$.
\eproof

Combining Theorems \ref{thm-1.8} and \ref{thm-1.18} we obtain
immediately that if $(X_n)$ is a sequence of Banach spaces with
uniformly bounded $K$-convexity constants and $\GL_2$-constants, then
$X=\left(\sum_{n=1}^\infty X_n\right)_2$ has $\GL_2$. However it was
pointed out to us by Junge that this conclusion can be obtained without
the $K$-convexity assumption by combining the inequality in
\ref{thm-1.17} with its dual form. We need

\begin{lemma}
\label{lemma-1.19}
Let $(X_n)$ be a sequence of Banach spaces, $X=\left(\sum_{n=1}^\infty
X_n\right)_2$, $P_n: X\to X_n$ the canonical projection.
\begin{itemize}
\item[(i)] If $T\in B(\ell_2,X)$ with $\sum\gamma_\infty(P_nT)^2<\infty$
then $T\in\Gamma_\infty(\ell_2,X)$ with
\[
\gamma_\infty(T)\le \left(\sum_{n=1}^\infty \gamma_\infty
(P_nT)^2\right)^{\frac12}.
\]
\item[(ii)] If $S\in B(X,\ell_2)$ with
$\sum_{n=1}^\infty\gamma_1(SP_n)^2<\infty$ then $S\in\Gamma_1(X,\ell_2)$
with
\[
\gamma_1(S) \le \left(\sum_{n=1}^\infty \gamma_1(SP_n)^2\right)^{\frac12}.
\]
\end{itemize}
\end{lemma}

\bproof
(i) follows immediately from Theorem \ref{thm-1.17} by applying trace
duality to the inequality there. Applying (i) to $X^*$ we obtain (ii).
\eproof

This leads to

\begin{theorem}
\label{thm-1.20}
Let $(X_n)$ be a sequence of Banach spaces all having $GL_2$ so that
$K=\sup_n \gl_2(X_n)<\infty$. Then $X=\left(\sum_{n=1}^\infty
X_n\right)_2$ has $\GL_2$.
\end{theorem}

\bproof
Let $T\in\Pi_1(X,\ell_2)$. From Theorem \ref{thm-1.17} and our
assumptions we get
\begin{equation}
\label{eq-1.21}
\sum_{n=1}^\infty \gamma_1(TP_n)^2\le K^2 \sum_{n=1}^\infty
\pi_1(TP_n)^2\le K^2\pi_1(T)^2.
\end{equation}
Lemma \ref{lemma-1.19} now gives that $T\in\Gamma_1(X,\ell_2)$ with
\begin{equation}
\label{eq-1.22}
\gamma_1(T) \le K\pi_1(T),
\end{equation}
which shows that $X$ has $\GL_2$.
\eproof

Let us end this section by discussing the following problem which seems
to be important since it has some applications to various areas of
Banach space theory.

\begin{problem}
\label{prob-1.21}
Let $(X_n)$ be a sequence of Banach spaces. Under which assumptions on
the $X_n$'s does there exist a constant $K$ so that
\begin{equation}
\label{eq-II}
\pi_1(T)\le K\left(\sum_{n=1}^\infty
\pi_1(TP_n)^2\right)^{\frac12}\quad\mbox{for all $T\in X\otimes\ell_2$}.
\end{equation}
\end{problem}

The next theorem gives some conditions for the inequality (\ref{eq-II})
to hold. (iii) was shown to us by Junge.

\begin{theorem}
\label{thm-1.22}
Let $(X_n)$ be a sequence of Banach spaces, $X=(\sum X_n)_2$. The
inequality (\ref{eq-II}) holds, if one of the following conditions is
satisfied.
\begin{itemize}
\item[(i)] $X^*_n$ has $\GAP$ for every $n\in\bN$ with uniformly bounded
$\GAP$-constants.
\item[(ii)] $X_n$ has $\GL_2$ for every $n\in\bN$ and $\sup
\gl_2(X_n)<\infty$.
\item[(iii)] $X_n$ is of cotype 2 for every $n\in \bN$ with uniformly
bounded cotype 2 constants.
\end{itemize}
\end{theorem}

\bproof
If (i) is satisfied, we choose $K\ge 0$ so that
\[
\ell(S) \le K\pi_1(S^*) \quad\mbox{for all $S\in\ell_2\otimes X^*_n$.}
\]
$X^*$ has $\GAP$ by Theorem \ref{thm-1.18} and by repeating the
calculations there with $X$ replaced by $X^*$ combined with
Proposition \ref{prop-0.3} we get for every $T\in X\otimes\ell_2$
\begin{equation}
\label{eq-1.23}
\pi_1(T)\le K_1\ell(T^*) = K_1\left(\sum_{n=1}^\infty
\ell(P^*_nT^*)^2\right)^{\frac12} \le KK_1\left(\sum_{n=1}^\infty
\pi_1(TP_n)^2\right)^{\frac12}
\end{equation}
which gives (\ref{eq-II}).

Assume next that (ii) holds. Put $K=\sup_n \gl_2(X_n)$. If $K_G$ denotes
the Grothendieck constant, then by repeating the calculations in the
proof of Theorem \ref{thm-1.20} we get for every $T\in X\otimes\ell_2$
\begin{equation}
\label{eq-1.24}
\pi_1(T) \le K_G\gamma_1(T)\le K_G\left(\sum_n
\gamma_1(TP_n)^2\right)^{\frac12} \le KK_G\left(\sum_{n=1}^\infty
\pi_1(TP_n)^2\right)^{\frac12}
\end{equation}
which gives (\ref{eq-II}).

Assume finally that $X$ is of cotype 2 with constant $K$ and let
$S\in\ell_2\otimes X$. By \cite[Theorem 12.2]{NT} we get that
\begin{equation}
\label{eq-1.25}
\pi_2(P_nS)\le K\ell(P_nS) \quad\mbox{for all $n\in\bN$}
\end{equation}
and hence
\begin{equation}
\label{eq-1.26}
\left(\sum_{n=1}^\infty \pi_2(P_nS)^2\right)^{\frac12} \le
K\left(\sum_{n=1}^\infty \ell(P_nS)^2\right)^{\frac12} = K\ell(S) \le
K_2K\pi_2(S),
\end{equation}
where $K_2$ is the constant from Proposition \ref{prop-0.3}.

Dualizing (\ref{eq-1.26}) and using again that $X$ is of cotype 2 we
obtain for every $T\in X\otimes\ell_2$
\begin{equation}
\label{eq-1.27}
\pi_1(T)\le K\pi_2(T)\le K^2K_1\left(\sum_{n=1}^\infty
\pi_2(TP_n)^2\right)^{\frac12} \le K^2K_1\left(\sum_{n=1}^\infty
\pi_1(TP_n)^2\right)^{\frac12},
\end{equation}
which gives (\ref{eq-II}).
\eproof

\setcounter{equation}{0}

\section{$\GAP$ and Extension Properties of Certain Classes of Operators}
\label{sec-2}

In this section we shall prove some results concerning extensions of
certain operators defined on a Banach space with $\GAP$ with values in a
Hilbert space. We start with the following

\begin{theorem}
\label{thm-2.1}
Let $X$ be a Banach space with $\GAP$. Then there is a constant $K$ so that
for every subspace $E\subseteq X$ and every $T\in\ell_2\otimes E$ we
have
\begin{equation}
\label{eq-2.1}
\pi_1(T^*)\le K\pi_1(T^*Q)
\end{equation}
where $Q$ is the canonical quotient map of $X^*$ onto $E^*$.

Consequently, every $S\in\Gamma_1(E,\ell_2)$ admits an extension
$\widetilde{S}\in\Gamma_1(X,\ell_2)$ with
\begin{equation}
\label{eq-2.2}
\gamma_1(\widetilde{S})\le K\gamma_1(S).
\end{equation}
\end{theorem}

\bproof
Let $C$ be the $\GAP$-constant of $X$ and let $T\in\ell_2\otimes E$ be
arbitrary. It is obvious that $\ell(T: \ell_2\to E)=\ell(T: \ell_2\to
X)$ and hence
\begin{equation}
\label{eq-2.3}
\pi_1(T^*) \le K_1\ell(T)\le K_1C\pi_1(T^*Q),
\end{equation}
where $K_1$ is the constant from Proposition \ref{prop-0.3},
(\ref{eq-2.3}) gives (\ref{eq-2.1}) with $K=K_1C$.

An easy dualization argument shows that the second statement is
equivalent to $\Gamma^*_1(\ell_2,E)\subseteq \Gamma_1^*(\ell_2,X)$ with
$K$-equivalence between the norms. ($\Gamma_1^*$ denotes the dual
operator ideal.)

However, $\Gamma^*_1(\ell_2,E)=\{T\in B(\ell_2,E)\mid
T^*\in\Pi_1(E^*,\ell_2)\}$ and similarly for $X$, and hence the latter
statement is exactly (\ref{eq-2.1}).
\eproof

The next theorem gives a characterization of subspaces $E$ of a given
Banach space $X$ so that $E^*$ has $\GAP$ in terms of extensions of
1-summing operators.

\begin{theorem}
\label{thm-2.2}
Let $X$ be a Banach space and $E$ a subspace. Consider the statements
\begin{itemize}
\item[(i)] $E^*$ has $\GAP$.
\item[(ii)] There exists a constant $K\ge 0$ so that every
$T\in\Pi_1(E,\ell_2)$ admits an extension
$\widetilde{T}\in\Pi_1(X,\ell_2)$ with $\pi_1(\widetilde{T})\le
K\pi_1(T)$.
\end{itemize}

If $X$ is of finite cotype then (i) implies (ii). If $X^*$ has $\GAP$ then
(ii) implies (i).
\end{theorem}

\bproof
By duality (ii) is equivalent to 

(iii) $\Gamma_\infty(\ell_2,E)\subseteq
\Gamma_\infty(\ell_2,X)$ with equivalence between the norms.

Let $X$ be of finite cotype and assume that $E^*$ has $\GAP$. 

We wish to show that (iii) holds. By Proposition \ref{prop-1.5} and Theorem
\ref{thm-1.6} there exist constants $K\ge 0$ and $C\ge 0$ so that if
$T\in\ell_2\otimes E$,
\begin{equation}
\label{eq-2.4}
\gamma_\infty(T)\le C\ell(T)\le KC\gamma_\infty(T: \ell_2\to X)
\end{equation}
which shows that (iii) holds.

Assume next that $X^*$ has $\GAP$ with constant $M$ and that (ii) holds. It
clearly follows that there is a constant $K\ge 0$ so that every $T\in
\Pi_1(E,\ell_2)$ admits an extension $\widetilde{T}\in\Pi_1(X,\ell_2)$
with
\begin{equation}
\label{eq-2.5}
\pi_1(\widetilde{T})\le K\pi_1(T).
\end{equation}
Let now $T=\sum_{j=1}^n f_j^*\otimes e_j\in E^*\otimes\ell_2$ and
let $\widetilde{T}$ be an extension of $T$ so that (\ref{eq-2.5})
holds. Without loss of generality we may assume that the range of
$\widetilde{T}$ is contained in $[e_j\mid 1\le j\le n]$ and since $X^*$
has $\GAP$ we therefore easily obtain
\begin{equation}
\label{eq-2.6}
\ell(T^*)\le \ell(\widetilde{T}^*)\le M\pi_1(\widetilde{T})\le
KM\pi_1(T),
\end{equation}
which shows that $E^*$ has $\GAP$.
\eproof

Combining Theorem \ref{thm-2.2} with the results of the previous section
we obtain the following corollary

\begin{corollary}
\label{cor-2.3}
Let $X$ be a Banach space of finite cotype with $\GL_2$ and let
$E\subseteq X$ be a subspace. Then the
following statements are equivalent.
\begin{itemize}
\item[(i)] $E$ has $\GL_2$.
\item[(ii)] Every operator $T\in\Pi_1(E,\ell_2)$ admits an extension $\widetilde{T}\in\Pi_1(X,\ell_2)$.
\end{itemize}
\end{corollary}

\bproof
Trivially (ii) implies(i) (for this the finite cotype assumption on $X$
is superfluous). Assume next that $E$ has $\GL_2$ and let $T\in
\Pi_1(E,\ell_2)$; hence $T\in\Gamma_1(E,\ell_2)$ as well and since $X$
has $\GAP$ we get from Theorem \ref{thm-2.1} that $T$ admits an extension
$\tilde{T}\in\Gamma_1(X,\ell_2)\subseteq\Pi_1(X,\ell_2)$.
\eproof

The assumption that $X$ is of finite cotype cannot be omitted in
Corollary \ref{cor-2.3} as the following example shows:

\begin{example}
\label{ex-2.4}
Let $E$ be a subspace of $\ell_\infty$ isometric to $\ell_1$, and let
$T\in B(E;\ell_2)$ be onto. $E$ has $GL_2$ and $T$ is absolutely summing
by Grothendieck's theorem. If $T$ could be extended to a
$\tilde{T}\in\pi_1(\ell_\infty,\ell_2)$, then $\ell(\tilde{T}^*)<\infty$
since $\ell_\infty^*$ has $\GAP$. Since $\tilde{T}^*$ is an isomorphism this
is a contradiction.
\end{example}



\noindent Department of Mathematics, University of Missouri, \\
Columbia, MO 65211, USA\\
pete\makeatletter{@}\makeatother casazza.math.missouri.edu
  \ \\
Department of Mathematics and Computer Science, Odense University,\\
Campusvej 55, DK-5230 Odense M, Denmark\\
njn\makeatletter{@}\makeatother imada.ou.dk

\end{document}